\documentclass[11pt]{article}
\usepackage{amssymb}
\usepackage{verbatim,cite}
\usepackage{latexsym}
\usepackage{amsmath}
\usepackage{eucal}

\newtheorem{sub}{\name}[section]
\newcommand{\bs}{
\begin{sub}}
\newcommand{\es}{
\end{sub}}
\newcommand{\bsl}[1]{
\begin{sub}\label{#1}}
\newcommand{\bth}[1]{\def\name{Theorem}
\begin{sub}\label{t:#1}}
\newcommand{\blemma}[1]{\def\name{Lemma}
\begin{sub}\label{l:#1}}
\newcommand{\bcor}[1]{\def\name{Corollary}
\begin{sub}\label{c:#1}}
\newcommand{\bdef}[1]{\def\name{Definition}
\begin{sub}\label{d:#1}}
\newcommand{\bprop}[1]{\def\name{Proposition}
\begin{sub}\label{p:#1}}
\newcommand{\brem}[1]{\def\name{Remark}
\begin{sub}\label{r:#1}}
\newcommand{\bex}[1]{\def\name{Example}
\begin{sub}\label{e:#1}}


\newcommand{\BA}{
\begin{array}}
\newcommand{\EA}{
\end{array}}

\newcommand{\BAN}{\renewcommand{\arraystretch}{1.2}
\setlength{\arraycolsep}{2pt}
\begin{array}}
\newcommand{\BAV}[2]{\renewcommand{\arraystretch}{#1}
\setlength{\arraycolsep}{#2}
\begin{array}}
\newcommand{\BSA}{
\begin{subarray}}
\newcommand{\ESA}{
\end{subarray}}
\newcommand{\BAL}{
\begin{aligned}}
\newcommand{\EAL}{
\end{aligned}}
\newcommand{\BALG}{
\begin{alignat}}
\newcommand{\EALG}{
\end{alignat}}
\newcommand{\BALGN}{
\begin{alignat*}}
\newcommand{\EALGN}{
\end{alignat*}}
\newcommand{\note}[1]{\textit{#1.}\hspace{2mm}}
\newcommand{\Proof}{\note{Proof}}
\newcommand{\qeda}{\hspace{10mm}\hfill $\square$}

\newcommand{\Remark}{\note{Remark}}

\newcommand{\abs}[1]{\left |#1\right |}
\newcommand{\norm}[1]{\left \|#1\right \|}

\def\angb<#1>{\langle #1 \rangle}

\newcommand{\myfrac}[2]{{\displaystyle \frac{#1}{#2} }}
\newcommand{\myint}[2]{{\displaystyle \int_{#1}^{#2}}}

\newcommand{\prt}{
\partial}

\newcommand{\ti}{\times}

\newcommand{\nind}{\noindent}

\def\ga{\alpha}            
             \def\ge{\epsilon}
                         
\def\gf{\phi}           
            \def\gl{\lambda}
                 \def\gp{\pi}
    \def\gr{\rho}        
\def\gs{\sigma}       \def\gt{\tau}
      \def\gw{\omega}
                
     \def\Gd{\Delta}      \def\Gf{\Phi}

   \def\CE{{\mathcal E}}

   \def\BBN {\mathbb N}    
   \def\BBR {\mathbb R}

\allowdisplaybreaks
\begin{document}
\title {\bf  Admissible initial growth for  diffusion equations with weakly superlinear absorption}
\author{{\bf\large Andrey Shishkov}
\hspace{2mm}\vspace{3mm}\\
{\it \normalsize  Institute of Applied Mathematics and Mechanics of NAS of Ukraine},\\
{\it\normalsize R. Luxemburg str. 74, 83114 Donetsk, Ukraine}\\
\vspace{3mm}\\
{\bf\large Laurent V\'eron}
\vspace{3mm}\\
{\it\normalsize Laboratoire de Math\'ematiques et Physique Th\'eorique, CNRS UMR 6083},
\\
{\it\normalsize Universit\'e Fran\c{c}ois-Rabelais,  37200 Tours,
France}}

\date{}
\maketitle


\noindent{\small {\bf Abstract} We study the admissible growth at infinity of initial data of positive solutions of $\prt_t u-\Gd u+f(u)=0$ in $\BBR_+\ti\BBR^N$ when $f(u)$ is a continuous  function, {\it mildly} superlinear at infinity, the model case being $f(u)=u\ln^\ga (1+u)$ with $1<\ga<2$. We prove in particular that if the growth of the initial data at infinity  is too strong, there is no more diffusion and the corresponding solution satisfies the ODE problem $\prt_t \gf+f(\gf)=0$ on $\BBR_+$ with $\gf(0)=\infty$.}
\tableofcontents

\section{Introduction and formulation of the results}
Let $h$ be a continuous nondecreasing function defined on $\BBR_+$ and vanishing only at $0$. It is well
known that for any continuous and bounded function $g$ belonging to  $C^+_b(\BBR^N)$, the cone of
bounded nonnegative continuous functions on $\BBR^N$, there exists a unique weak solution  $u:=u_g\in C^+_b(\BBR_+\times\BBR^N)$ of
 \begin{equation}\label{A1}\BA {lll}
\prt_tu-\Gd u +uh(u)=0\qquad&\text{in }\,Q_{\BBR^N}^\infty:=\BBR_+\times\BBR^N,\\[2mm]\phantom{--;--}
\!\displaystyle\lim_{t\to 0}u(t,.)=g\qquad&\text{locally uniformly in }\BBR^N.
\EA\end{equation}
Furthermore, the solution $u$ satisfies
 \begin{equation}\label{A2}\BA {lll}
0\leq u(x,t)\leq \Gf_{\norm {g}_{L^\infty}}(t)\qquad\forall
(t,x)\in Q_{\BBR^N}^\infty, \EA\end{equation} where $\Gf_a$ is the solution of the following Cauchy problem:
 \begin{equation}\label{A3}\BA {lll}\phantom{--}
\Gf_t+\Gf h(\Gf)=0\qquad\text{on }\BBR_+,\\
\phantom{_t\Gf+-\Gf-}
\Gf(0)=a.
\EA\end{equation}
Since  $h(0)=0$ it follows easily that
$$\Gf_a(t)\geqslant 0\quad \forall t>0,\forall a\geqslant 0.$$

Moreover the family $\{\Gf_a(t)\}$ is monotonically increasing with respect to the parameter $a$, and the condition:
\begin{equation}\label{H1}\BA {l}
\myint{c}{\infty}\myfrac{ds}{sh(s)}<\infty,\quad c=\textrm{const} >0,
\EA\end{equation}
is equivalent to the existence  of a solution $\Gf_{\infty}(t)$   of equation in $(\ref{A3})$  such that $\Gf_{\infty}(t)<\infty$, $\forall t>0$ and with infinite initial data, i. e. 
$$\lim_{t\rightarrow 0}\Gf_{\infty}(t)=\infty.$$

Now we come to our main subject, to study problem $(\ref{A1})$  with an  initial data $g(x)$ unbounded and tending to infinity at infinity. It is clear that the character of growth of $h(s)$ at infinity defines the class of initial functions $g$ {\it of solvability} of problem under consideration. For example, if $h(s)$ is bounded, then the corresponding class of solvability is the Tikhonov class \cite{Tikh} $\{g:g(x)\leqslant c \exp(c_1|x|^2),c,c_1=\textrm{const}<\infty\}$. When $h(s)$ tends to infinity at infinity, the class of admissible initial data is larger that the  Tikhonov class. 
If $h(s)$ increases at infinity fast enough in the sense that condition $(\ref{H1})$  is satisfied, then problem  $(\ref{A1})$   is solvable for any nonnegative continuous function $g$ in the sense that there exists a {\it prospective minimal solution} $\underline u_g$ which is the limit when $n\to\infty$ of the solutions $u_{g_n}$ of
 \begin{equation}\label{A4}\BA {lll}\phantom{-}
\prt_tu-\Gd u +uh(u)=0\qquad&\text{in }\,Q_{\BBR^N}^\infty\\[2mm]\phantom{--..--}
\!\displaystyle\lim_{t\to
0}u(t,.)=g\chi_{_{B_n}}\qquad&\text{in }L^1(\BBR^N),
\EA\end{equation} where $B_n$ denotes the open ball of radius $n$
and $\chi_{A}$ is characteristic function of the set $A$,
and there holds
 \begin{equation}\label{A5--}
 0\leq \underline u_g\leq \Gf_\infty.
 \end{equation}
{\it But the main question is to know whether the prospective minimal solution is truly a solution with initial data  $g(.)$}. If it is the case we say that the prospective minimal solution is the minimal solution. Another assumption on $h$ which plays a fundamental role in the study is the so-called Keller-Osserman condition (see \cite{Ke}, \cite{Oss}),
\begin{equation}\label{H1-1}\BA {l}
\myint{a}{\infty}\myfrac{ds}{\sqrt{H(s)}}<\infty\;\;\text{for some }\,a>0,\,\text{ where }\; H(t)=\myint{0}{s}th(t)dt.
\EA\end{equation}
When this condition is satisfied, for any $R>0$ there exist solutions of 
\begin{equation}\label{H1-2}\BA {l}
-\Gd u+uh(u)=0\qquad\text{in }\; B_R\\
\phantom{--}
\displaystyle \lim_{|x|\to R}u(x)=\infty.
\EA\end{equation}

This condition gives rise to a localization phenomenon thanks to which we  prove an existence and uniqueness result of the solution with initial data $g$.\smallskip

 \nind{\bf Theorem A} {\it Assume $r\mapsto rh(r)$ is convex and satisfies
$(\ref{H1-1})$.
Then for any $g\in C^+(\BBR^N)$, $\underline u_g$ is the minimal solution with initial data $g$. Furthermore it is the unique nonnegative solution of problem $(\ref{A1})$}.\smallskip

In the class of functions $h(s)$ of the form $h(s)=s \ln^{\alpha}(1+s)$, condition $(\ref{H1})$  is equivalent to $\alpha>1$, but condition $(\ref{H1-1})$  is equivalent to $\alpha>2$. In general it is easy to show that condition $(\ref{H1-1})$  is stronger than condition $(\ref{H1})$ .

When $h(s)$ is a power function, the class of existence and uniqueness is much larger than the class of Th. A. A complete description of this existence and uniqueness class is based upon the notion of initial trace which has been thoroughly investigated by Marcus and V\'eron \cite{MaVe1}, \cite{MaVe2} and Gkikas and V\'eron \cite{GkVe}.\\

When
\begin{equation}\label{H2}\BA {l}
\myint{a}{\infty}\myfrac{ds}{\sqrt{H(s)}}=\infty\qquad\forall a>0,
\EA\end{equation} uniqueness may not hold in the class of
unbounded solutions. If, for any $b>0$,  $V_b$ denotes the maximal
solution of the following Cauchy problem
 \begin{equation}\label{A5-}\BA {lll}\phantom{-}
V_{rr}+\myfrac{N-1}{r}V_{r}-Vh(V)=0\qquad\text{on }\,(0,R_b)\\
\phantom{_{rr}+\myfrac{N-1}{r}_{r}V,,Vh(V)}
\!\!V(0)=b\\
\phantom{_{rr}+\myfrac{N-1}{r}_{r}VVh(V)}
V_{r}(0)=0,
\EA\end{equation}
then $R_b=\infty$. Actually, multiplying $(\ref{A5-})$  by $V_r$, we get easily
$$2^{-1}\frac{d}{dr}|V_r|^2=\frac{d}{dr}H(V)-\frac{N-1}{r}|V_r|^2\leqslant\frac{d}{dr}H(V).$$
Since $V_r(0)=0$, we derive
$$V_r(r)\leqslant\sqrt{2}\sqrt{H(V(r))}\qquad \forall r>0.$$
Integrating this last inequality we obtain the {\it a priori} estimate
$$V(R)=V_b(R)\leqslant\bar{V_b}(R)\qquad \forall R>0,$$
where the function $\bar{V_b}(R)$ is defined by the identity:
$$F_b(\bar{V_b}(R))=\sqrt{2}R \qquad \forall R>0,\text{ where }\; F_b(v):=\int_b^v\frac{ds}{\sqrt{H(s)}},$$
 (see e.g. \cite{VaVe} also). Moreover it is easy to see that for arbitrary $a>b>0$, $V_a(r)\geqslant V_b(r)$ $\forall r>0$. Actually, due to the monotonicity of $h$, there holds
$$V_{a_{rr}}(0)=\frac{1}{N}V_a(0)h(V_a(0))=\frac{1}{N}ah(a)>\frac{1}{N}bh(b)=V_{b_{rr}}(0),$$
from  $(\ref{A5-})$.
Since $V_{a_{r}}(0)=V_{b_r}(0)=0$ it follows from this last inequality that the function  $r\mapsto W(r)= V_a(r)-V_b(r)$
is increasing near $r=0$; it remains increasing on whole $\BBR_+$, since, if we assume that there exists
 $r_0>0$ where $W$ reaches a local maximum, then $W_r(r_0)=0,W_{rr}(r_0)\leq 0$, but from equation $(\ref{A5-})$  we have:
$$W_{rr}(r_0)=V_a(r_0)h(V_a(r_0))-V_b(r_0)h(V_b(r_0))>0,$$
which is a contradiction. Furthermore, Nguyen Phuoc and V\'eron proved in \cite{NgVe} that if $g$ satisfies
\begin{equation}\label{A5}\BA {l}
V_c(|x|)\leq g(x)\leq V_b(|x|)\qquad\forall x\in\BBR^N,
\EA\end{equation}
for some $b>c>0$, then there exists at least two different solutions of $(\ref{A1})$ defined in $Q_{\BBR^N}^\infty$: the minimal one $\underline u_g$ which satisfies
 \begin{equation}\label{A6}\BA {l}
 \underline u_g(x,t)\leq \Phi_\infty(t)\qquad\forall (x,t)\in Q_{\BBR^N}^\infty,
\EA\end{equation}
and another one $u_g$ such that
 \begin{equation}\label{A7}\BA {l}
V_c(|x|)\leq u_g(x,t)\leq V_b(|x|)\qquad\forall (x,t)\in Q_{\BBR^N}^\infty.
\EA\end{equation}

It is not clear wether there exists a maximal solution or
not. However, if $g$ satisfies $(\ref{A5})$, then there exists a 
minimal solution $\underline u_{g,c,b}$ and a maximal one
$\overline u_{g,c,b}$ in the class $\CE_{c,b}(g)$ of solutions
 of problem $(\ref{A1})$, satisfying inequalities $(\ref{A7})$.
 These two solutions can be constructed by the following
 approximate scheme: we define the sequence $\{\underline u_n\}$
 of solutions of the Cauchy-Dirichlet problem
 \begin{equation}\label{A8-1}\BA {lll}
\prt_t u-\Gd  u + uh( u)=0\qquad&\text{in }\,Q_{B_n}^\infty:=\BBR_+\times B_n\\[2mm]\phantom{------}
\!\! u(t,x)=V_c(n)\qquad&\text{in }\,\prt_\ell Q_{B_n}^\infty:=\BBR_+\times \prt B_n\\[2mm]\phantom{------}
\!\displaystyle
 u(0,.)=g\qquad&\text{in }B_n;
\EA\end{equation} 
then it is easy to check using comparison
principle that the sequence $\{\underline u_n\}$ is increasing and
converges to $\underline u_{g,c,b}$. Similarly, the sequence
$\{\overline u_n\}$ of solutions of the same equation in
$Q_{B_n}^\infty$ with the same initial data and boundary value
$V_b(n)$ is decreasing and converges to $\overline
u_{g,c,b}$.\smallskip

In this paper we consider the case where the initial data $g$ grows at infinity faster than any function $V_b$ with arbitrary $b<\infty$. Our aim is to describe analogs of the "maximal" solution
$u_g$ from $(\ref{A7})$  and prospective minimal solution $\{\underline u_g\}$ from  $(\ref{A4})$, $(\ref{A5--})$. For any $a>0$ we denote by $u:=u_{a,n}$ the solution of
 \begin{equation}\label{A8}\BA {lll}
\prt_tu-\Gd  u + uh( u)=0\qquad&\text{in }\,Q_{B_n}^\infty:=\BBR_+\times B_n\\[2mm]\phantom{------}
\!\!u(t,x)=V_a(n)\qquad&\text{in }\,\prt_\ell Q_{B_n}^\infty:=\BBR_+\times \prt B_n\\[2mm]\phantom{------}
\!\displaystyle
 u(0,.)=\min\{V_a,g\}\qquad&\text{in }B_n,
\EA\end{equation}

\smallskip

Due to the comparison principle it is clear that  $u_{a,n}\leqslant V_a$ in $Q_{B_n}^\infty$. The next result highlights a phenomenon of 
{\it instantaneous blow-up} of the maximal solution if the initial data grows too fast at infinity.\smallskip

 \nind{\bf Theorem B} {\it Assume $r\mapsto rh(r)$ is convex and satisfies $(\ref{H1})$ and $(\ref{H2})$. If  $g\in C^+(\BBR^N)$,  satisfies
\begin{equation}\label{A9}\BA {l}
\displaystyle\lim_{|x|\to\infty}\myfrac{g(x)}{V_a(|x|)}=\infty\qquad\forall a>0,
\EA
\end{equation}
then for arbitrary $m\in\mathbb{N}$ the sequence $\{u_{a,n}\}_{n>m}$ decreases and converges in $Q_{B_m}^\infty$ to a function $u_a$ which is solution of $(\ref{A1})$ with initial data $\min\{V_a,g\}$. Furthermore
$u_a(t,x)\to\infty$ for any $(t,x)\in Q_{\BBR^N}^\infty$ as $a\to\infty$. Thus, the function identically  equal to $\infty$ in $Q^\infty_{\BBR^N}$ can
be considered as the "maximal" solution of problem $(\ref{A1})$ in the case  of  $(\ref{A9})$}.
\medskip

Let us remark that in subsection 3.1 we find the asymptotic expression of the functions $V_a$ for the model nonlinearities $h$,
which makes the condition  $(\ref{A9})$  more explicit.\smallskip

A fundamental example of equations with nonlinearities satisfying $(\ref{H1})$ and $(\ref{H2})$ is provided by
 \begin{equation}\label{A10}\BA {lll}
\prt_tu-\Gd  u + u\ln^\ga(1+u)=0\qquad&\text{in }\,Q_{\BBR^N}^\infty
\EA\end{equation}
with $1<\ga\leq 2$. With this specific type of nonlinearity we prove:

\smallskip

 \nind{\bf Theorem C} {\it Assume $1<\ga< 2$ and  $g\in C^+(\BBR^N)$,  satisfies condition  $(\ref{A9})$,
which due to Proposition 3.1 has now the following form
$$\lim_{|x|\rightarrow\infty}g(x)\exp\left(-c_{\alpha}|x|^{\frac{2}{2-\alpha}}\right)=\infty,\quad c_{\alpha}=\left(\frac{2-\alpha}{2}\right)^{\frac{2}{2-\alpha}}.$$
Then the prospective minimal solution $\underline u_g$ of $(\ref{A10})$ with initial data $g$ is $\Gf_\infty$.}\smallskip

Notice that the two types of generalized approximative solutions of problem  $(\ref{A1})$, obtained in Theorems B, C, "forget"
the real initial condition from    $(\ref{A1})$: in another words, they realize infinite initial jump.

\section{The maximal solution}
\setcounter{equation}{0}
\subsection{Proof of Theorem A }

The fact that $\underline u_g$ is a solution of $(\ref{A1})$, and clearly the minimal one, is more or less standard, but we recall its proof for the sake of completeness since it contains the localization principle. For $m\in\BBN^*$ let $v_m$ be the minimal  solution of
\begin{equation}\label{B1}\BA {l}
-\Gd v+vh(v)=0\qquad\text{in }\;B_m\\\phantom{--}
\displaystyle \lim_{|x|\to m}\!v(x)=\infty.
\EA
\end{equation}
Such a solution exists by \cite{Ke} or \cite{Oss} because
$(\ref{H1-1})$ holds. It is nonnegative and radial as limit of the nonnegative radial
functions $v_{m,k}$, $k\in\BBN^*$ which are the solutions of
$(\ref{B1})$ with finite boundary data $v_{m,k}=k$ on $\prt B_m$.
Moreover  $v_{m,k}$, and thus $v_{m}$, is an increasing function
of $|x|$. Clearly $v_m\geq 0$ and it is a stationary solution of
$(\ref{A1})$ in $Q^\infty_{B_m}$. For $n\geq m$, let $u_{g_n}$  be the solution of $(\ref{A4})$ and $\gamma_m=\max\{g(x):|x|\leq m\}$. Then $v_m+\gamma_m$ is a super solution of $(\ref{A4})$ in $Q_{B_m}^\infty$ which dominates $u_n$ on $\prt_\ell Q_{B_m}^\infty\cup \{0\}\times B_m$ Thus $v_m+\gamma_m\geq u_n$ in $Q_{B_m}^\infty$. The set of functions $\{u_n\}$ is bounded  and uniformly continuous in $Q^T_{B_{m-1}}$, by standard regularity theory for parabolic equations,  thus it converges  uniformly therein to $\underline u_g$ and $\underline u_g\lfloor_{B_m\ti\{0\} }=g$. This implies that $\underline u_g$ has $g$ as initial data. \smallskip

Assume now that $u$ another solution with the same initial data $g$. We set $w=u-\underline u_g$. Since $r\mapsto rh(r)$ is convex and $u-\underline u_g$ is positive,
$$uh(u)\geq \underline u_gh(\underline u_g)+(u-\underline u_g)h(u-\underline u_g).$$
Therefore $w$ is a subsolution of problem $(\ref{A1})$, and $w(t,x)\to 0$ as $t\to 0$, locally uniformly in $
\BBR^N$.  By the comparison principle
$$
w(t,x)\leq v_n(x)\qquad\text{in }\;Q^\infty_{B_n},
$$
where $v_n$ satisfies $(\ref{B1})$ in $B_n$. Furthermore $n\mapsto v_n$ is decreasing with limit $v_\infty$ as $n\to\infty$. The function $v_\infty$ verifies
$$-\Gd v+vh(v)=0\qquad\text{in }\;\BBR^N.
$$
Furthermore it is nonnegative, radial and nondecreasing with respect to  $|x|$. In order to prove that $v=0$, we return to $v_{n}$ which satisfies
$$
v_{n\,r}=r^{1-N}\myint{0}{r}s^{N-1}v_n(s)h(v_n(s))ds\leq
v_n(r)h(v_n(r))r^{1-N}\myint{0}{r}s^{N-1}ds=\myfrac{r}{N}v_n(r)h(v_n(r)).
$$
Thus
$$-v_{n\,rr}+v_n(r)h(v_n(r))=\myfrac{N-1}{r}v_{n\,r}\leq \left(1-\myfrac{1}{N}\right)v_n(r)h(v_n(r))
$$
which implies
$$-v_{n\,rr}+\myfrac{1}{N}v_n(r)h(v_n(r))\leq 0.
$$
Integrating twice yields
\begin{equation}\label{B}
\myint{v_n(r)}{\infty}\myfrac{ds}{\sqrt{H(t)}}\geq \sqrt\myfrac{2}{N}(n-r),
\end{equation}
where $H$ has been defined in $(\ref{H1-1})$. If we had $v_{\infty}(r)>0$ for some $r>0$, it would imply
$$\infty>\myint{v_{\infty}(r)}{\infty}\myfrac{ds}{\sqrt{2H(t)}}\geq \infty,
$$
a contradiction. Thus $v_\infty(r)=0$ and $w(t,x)=0$.\qeda \medskip

\subsection{Proof of Theorem B }

We recall that $(\ref{A9})$ holds and that $u_{a,n}$ denotes the solution of $(\ref{A8})$. Since $V_a\lfloor_{Q^\infty_{B_n}}$ is the solution of the Cauchy-Dirichlet problem
 \begin{equation}\label{B2}\BA {lll}
\prt_tu-\Gd  u + uh( u)=0\qquad&\text{in }\,Q_{B_n}^\infty:=\BBR_+\times B_n\\[2mm]\phantom{------}
\!\!u(t,x)=V_a(n)\qquad&\text{in }\,\prt_\ell Q_{B_n}^\infty:=\BBR_+\times \prt B_n\\[2mm]\phantom{------}
\!\displaystyle
 u(0,.)=V_a\qquad&\text{in }B_n,
\EA\end{equation}
it is larger than $u_{a,n}$. Thus $u_{a,n+1}\lfloor_{\prt_\ell {Q^\infty_{B_n}}}\leq u_{a,n}\lfloor_{\prt_\ell {Q^\infty_{B_n}}}=V_a$. Since $u_{a,n}(0,.)=u_{a,n+1}\lfloor_{B_n}(0,.)$ it follows that $u_{a,n+1}\lfloor_{ {Q^\infty_{B_n}}}\leq u_{a,n}$. Then  $\{u_{a,n}\}$ is a decreasing sequence, and its limit $u_a$ is a solution of $(\ref{A1})$, which the first claim. By the same argument, $u_{a,n}\leq u_{b,n+1}\lfloor_{ Q^\infty_{B_n}}$ in $Q^\infty_{B_n}$ for $b>a$. Hence $u_a\leq u_b$.
We introduce the sequence $\{r_a\}:r_a\to \infty$ as $a\to
\infty$ defined by:
\begin{equation}\label{B3}
r_a=\inf\{r>0:g(x)\geqslant V_a(x)\quad \forall\,|x|\geqslant r\},
\end{equation}
and, for $n\geq r_a$, we set $w_{a,n}=V_a-u_{a,n}$.
By convexity $w_{a,n}$ satisfies
 \begin{equation}\label{B4}\BA {lll}
\prt_tw_{a,n}-\Gd  w_{a,n} + w_{a,n}h( w_{a,n})\leq
0\qquad&\text{in }\,Q_{B_n}^\infty:=\BBR_+\times B_n,
\\[2mm]\phantom{---------,}
\!\!w_{a,n}(t,x)=0\qquad&\text{in }\,\prt_\ell Q_{B_n}^\infty:=\BBR_+\times \prt
B_n,\\[2mm]\phantom{---------,}
\!\displaystyle w_{a,n}(0,x)=(V_a-g)_+\qquad&\text{in }B_n.
\EA\end{equation}Therefore
 \begin{equation}\label{B5}
w_{a,n}(t,x)<\Gf_\infty(t) \qquad\text{in }\,Q_{B_n}^\infty,
\end{equation}
where $\Gf_\infty$ is defined in $(\ref{A3})$ with $a=\infty$. Actually,
 \begin{equation}\label{B6}
\myint{\Gf_\infty(t)}{\infty}\myfrac{ds}{sh(s)}=t.
\end{equation}
Notice also that the sequence $\{w_{a,n}\}$ is increasing and it converges, as  $n\to\infty$, to $w_a=V_a-u_a$, which is dominated by $\Gf_\infty$
Thus
 \begin{equation}\label{B7}
u_a(t,x)\geq V_a(x)-\Gf_\infty(t)\geq a-\Gf_\infty(t) \qquad\text{in }\,Q_{\BBR^N}^\infty.
\end{equation}
Letting $a\to\infty$ implies the claim.\qeda

\section{The prospective minimal solution}
\setcounter{equation}{0}
In this section we consider equation $(\ref{A10})$ with $1<\ga< 2$.
\subsection{The stationary problem}
\bprop{Stat} Assume $1<\ga < 2$, $a>0$ and $V_a$ is the solution
of
 \begin{equation}\label{m1}\BA {lll}
V_{rr}+\myfrac{N-1}{r}V_r-V\ln^\ga(V+1)=0\qquad\text{in }\BBR_+\\
\phantom{,\myfrac{N-1}{r}V_r-V\ln^\ga(V+1)}
V_r(0)=0
\\
\phantom{+\myfrac{N-1}{r}V_r-V\ln^\ga(V+1)}
V(0)=a.
\EA\end{equation}
Then
 \begin{equation}\label{m2}\BA {lll}
V_a(r)=e^{c_\ga r^{\frac{2}{2-\ga}}+O(1)}\qquad\text {as }\,r\to\infty,
\EA\end{equation}
where $c_\ga=\left(\frac{2-\ga}{2}\right)^{\frac{2}{2-\ga}}$.
\es
\Proof We write $W=\ln(V+1)$. Since $V$ is increasing $W>0,W_r\geq 0$ and
 \begin{equation}\label{m3}\BA {lll}
W_{rr}+W^2_r+\myfrac{N-1}{r}W_r-(1-e^{-W})W^\ga=0\qquad\text{in }\BBR_+.
\EA\end{equation}
Thus
$$W_{rr}+W^2_r-(1-e^{-W})W^\ga\leq 0.
$$
If we set $\gr=W$ and $ p(\gr)=W_r(r)$, then $\gr\in [a,\infty)$ and
$$pp'+p^2-(1-e^{-\gr})\gr^\ga\leq 0.
$$
This is a linear differential inequality in the unknown $p^2$. Integrating yields
 \begin{equation}\label{m4}p^2(\gr)\leq 2e^{-2\gr}\myint{a}{\gr}(e^{2s}-e^{s})s^\ga ds=\gr^\ga+O(1).
\end{equation}
Thus $W_r(r)\leq W^{\frac{\ga}{2}}(r)+O(1)$ as $r\to\infty$ which implies
 \begin{equation}\label{m5}
 W(r)\leq c_\ga r^{\frac{2}{2-\ga}}+O(1)\qquad\text{as }\,r\to\infty.
\end{equation}
Due to $(\ref{m5})$ relation $(\ref{m4})$ yields also the
following inequality
$$
0<W_r\leq c_\ga^{\frac{\ga}{2}} r^{\frac{\ga}{2-\ga}}(1+o(1)).
$$
Since $W(r)\to \infty$ as $r\to \infty$, it follows from
$(\ref{m3})$  and $(\ref{m4})$ that for any
$\epsilon>0$  there exists $r_\ge>0$ such that
$$
W_{rr}+W^2_r\geq (1-\ge)W^\ga\qquad\text{on }\,[r_\ge,\infty).
$$
Integrating this ordinary differential inequality
we get
 \begin{equation}\label{m6}W(r)\geq (1-\ge))c_\ga r^{\frac{2}{2-\ga}}(1+o(1))\qquad\text{as }\,r\to\infty.
\end{equation}
Since $\ge$ is arbitrary, we derive
 \begin{equation}\label{m7}W(r)= c_\ga r^{\frac{2}{2-\ga}}(1+o(1))\qquad\text{as }\,r\to\infty.
\end{equation}
From the above estimates, we can improve  $(\ref{m6})$. Using
$(\ref{m4})$ and $(\ref{m7})$ we deduce from $(\ref{m3})$:
$$
pp'+p^2= (1-e^{-\gr})\gr^\ga-\myfrac{N-1}{r}W_r\geq (1-e^{-\gr})\gr^\ga -c\gr^{\ga-1},
$$
from which it follows easily
 \begin{equation}\label{m8}p^2(\gr)\geq 2e^{-2\gr}\myint{a}{\gr}e^{2s}(s^\ga-c's^{\ga-1}) ds=\gr^\ga+O(1),
\end{equation}
by l'Hospital rule. Combined with $(\ref{m7})$ and  $(\ref{m5})$, it
implies
\begin{equation}\label{m9}W(r)= c_\ga r^{\frac{2}{2-\ga}}+O(1)\qquad\text{as }\,r\to\infty.
\end{equation}
Returning to $V_a$, we derive
 \begin{equation}\label{m10}
 V_a(r)= e^{c_\ga r^{\frac{2}{2-\ga}}+O(1)}\qquad\text{as }\,r\to\infty.
\end{equation}\qeda\medskip

\nind\Remark If $\ga=2$, the same method yields
 \begin{equation}\label{m11}
 V_a(r)= e^{e^r+O(1)}\qquad\text{as }\,r\to\infty.
\end{equation}
\subsection{Proof of Theorem C}
We recall that the prospective minimal solution $\underline u_ g$ is the limit, when $n\to\infty$ of the (increasing) sequence of solutions $\{u_{g_{\ell_n}}\}$ of
 \begin{equation}\label{m12}\BA {lll}
\prt_tu-\Gd u+u\ln^\ga(u+1)=0\qquad&\text{in }Q^\infty_{\BBR^N}\\[2mm]
\phantom{\prt-\Gd v+v\ln^\ga u1}
u(0,.)=g\chi_{B_{\ell_n}}\qquad&\text{in }\BBR^N,
\EA\end{equation} where $\{\ell_n\}$ is any increasing sequence
converging to $\infty$. Furthermore, if we replace  $g$ by its
maximal radial minorant defined by $\tilde g(r):= \min_{\abs
x=r}g(x)$, it satisfies also $(\ref{A9})$. Because of $(\ref{A9})$
there exists a sequence $\{r_n\}$ tending to infinity such that
$$
r_n=\inf \{r>0:\tilde g(s)\geq V_n(s)\,\,\,\forall s\geq r\},
$$
then $\tilde g (r_n)=V_n(r_n)$.

\nind{\it Step 1: Estimate from below}. Put
$$
g_n(|x|)=\left\{\BA {ll}\min\{\tilde g(r_n),\tilde g(|x|)\}\qquad&\text{if }|x|< r_n\\
\tilde g(r_n)\qquad&\text{if }|x|\geq r_n.
\EA\right.
$$
Let $\underline u_{g_n}$ be the minimal solution of
 \begin{equation}\label{m13}\BA {lll}
\prt_tu-\Gd u+u\ln^\ga(u+1)=0\qquad&\text{in }Q^\infty_{\BBR^N}\\[2mm]
\phantom{\prt-\Gd u+u\ln^\ga u1}
u(0,.)=g_n\qquad&\text{in }\BBR^N.
\EA\end{equation}
Then $\underline u_{g_n}\leq \Gf_\infty$. For any sequence $\{\ell_k\}$ converging to infinity and any fixed $k$, there exists $n_k$ such that for $n\geq n_k$, there holds $g\chi_{B_{\ell_k}}\leq g_n$. Since the sequence $\{\underline u_{g_n}\}$ is increasing, its limit  $u_\infty$ is a solution of $(\ref{A3})$ in $Q_\infty^{\BBR^N}$ which is larger than $u_{g_{\ell_k}}$ for any $\ell_k$,  and therefore larger also than $\underline u_{\tilde g}$. However, since $g_n\leq \tilde g$, $u_\infty\leq \underline u_{\tilde g}$. This implies
 \begin{equation}\label{m13*}\BA {lll}
u_\infty= \underline u_{\tilde g}\leqslant \Phi_\infty.
\EA\end{equation}

Next, since $\underline u_{g_n}(0,x)\leq g(r_n)$, it follows that $\underline u_{g_n}(t,x)\leq g(r_n)$. Let $\gw_n=\Gf_{g(r_n)}$, i.e. the solution of  $(\ref{A3})$ with $a=g(r_n)$, then $\gw_n$ satisfies
$$
\myint{\gw_n(t)}{g(r_n)}\myfrac{ds}{sh(s)}=t,
$$
and $\underline u_{g_n}\geq w_n$ where $w_n$ is the minimal solution of
 \begin{equation}\label{m14}\BA {lll}
\prt_tw-\Gd w+w\ln^\ga(\gw_n+1)=0\qquad&\text{in }Q^\infty_{\BBR^N}\\[2mm]
\phantom{;\prt-\Gd v+v\ln^\ga \gw_n1}
w(0,.)=g_n\qquad&\text{in }\BBR^N.
\EA\end{equation}
If we set $w_n(t,x)=e^{-\int_0^t\ln^\ga(\gw_n(s)+1)ds}z_n(t,x)$, then
 \begin{equation}\label{m15}\BA {lll}
\prt_tz_n-\Gd z_n=0\qquad&\text{in }Q^\infty_{\BBR^N}\\[2mm]
\phantom{\prt_tz_n}
z_n(0,.)=g_n\qquad&\text{in }\BBR^N.
\EA\end{equation}
Since
$$\BA {lll}
z_n(t,x)=\myfrac{1}{(4\gp t)^{\frac{N}{2}}}\myint{\BBR^N}{}e^{-\frac{|x-y|^2}{4t}}g_n(y)dy,
\EA$$
we can write $w_n(t,x)=I_n(t,x)+J_n(t,x)$ where
 \begin{equation}\label{m17}\BA {lll}
I_n(t,x)=\myfrac{e^{-\int_0^t\ln^\ga(\gw_n(s)+1)ds}}{(4\gp t)^{\frac{N}{2}}}\myint{\abs y\leq r_n}{}e^{-\frac{|x-y|^2}{4t}}g_n(y)dy,
\EA\end{equation}
and
 \begin{equation}\label{m18}\BA {lll}
J_n(t,x)=\myfrac{e^{-\int_0^t\ln^\ga(\gw_n(s)+1)ds}\tilde g(r_n)}{(4\gp t)^{\frac{N}{2}}}\myint{\abs y> r_n}{}e^{-\frac{|x-y|^2}{4t}}dy.
\EA\end{equation}
Clearly
 \begin{equation}\label{m19}\BA {lll}J_n(t,x)\geq \myfrac{e^{-\int_0^t\ln^\ga(\gw_n(s)+1)ds}\tilde g(r_n)}{(4\gp t)^{\frac{N}{2}}}
\myint{\abs y> r_n+|x|}{}e^{-\frac{|y|^2}{4t}}dy
\\\phantom{J_n(t,x)}\geq
\myfrac{e^{-\int_0^t\ln^\ga(\gw_n(s)+1)ds}\tilde g(r_n)}{(4\gp
t)^{\frac{N}{2}}} \left(\myint{\abs z>
r_n+|x|}{}e^{-\frac{z^2}{4t}}dz\right)^N. \EA\end{equation} This
integral term can be estimated by introducing Gauss error function
 \begin{equation}\label{m20}\BA {lll}
\text{ercf}(x)=\myfrac{2}{\sqrt\gp}\myint{x}{\infty}e^{-z^2}dz.
\EA\end{equation}
In dimension $N$, it implies easily
 \begin{equation}\label{m20'}\BA {lll}J_n(t,x)\geq e^{-\int_0^t\ln^\ga(\gw_n(s)+1)ds}\tilde g(r_n)
 \left(\text{ercf}\left(\myfrac{r_n+|x|}{2\sqrt t}\right)\right)^N.
 \EA\end{equation}
 Since
 $$\text{ercf}(x)=\myfrac{e^{-x^2}}{x\sqrt t}(1+O(x^{-2}))\quad\text{as }\,x\to\infty,
 $$
 we derive
  \begin{equation}\label{m21}\BA {lll}J_n(t,x)\geq\myfrac{\tilde g(r_n)}{((r_n+|x|)^2t)^{\frac{N}{2}}}
e^{-\int_0^t\ln^\ga(\gw_n(s)+1)ds-\frac{N(r_n+|x|)^2}{4t}}\left(1+O\left(\frac{t}{r^2_n}\right)\right).
 \EA\end{equation}
We write $\tilde g(r)=\exp(\gamma (r))-1$ and set
$$A_n(t,x)=\gamma (r_n)-\int_0^t\ln^\ga(\gw_n(s)+1)ds-\frac{N(r_n+|x|)^2}{4t}-N\ln(r_n+|x|)-
\myfrac{N}{2}\ln t.
$$
In order to have an estimate on $\gw_n(s)$, we fix $t\leq 1$ and $\tilde g(r_n)\geq 1$. There exists $a_0\geq 1$ such that
$$\min\left\{\myfrac{\gw_a(t)}{\gw_a(t)+1}:0\leq t\leq 1,\,a\geq a_0\right\}\geq\myfrac{1}{2}.
$$
In such a range of $a$ and $t$,
$$\BA {ll}\gw'+\gw\ln^\ga(\gw+1)=\gw'+\myfrac{\gw}{\gw+1}(\gw+1)\ln^\ga(\gw+1)\\[2mm]
\phantom{\gw'+\gw\ln^\ga(\gw+1)}
\geq \gw'+\myfrac{1}{2}(\gw+1)\ln^\ga(\gw+1),
\EA$$
which yields
 $$\ln^\ga(\gw_n(s)+1)\leq \left(\myfrac{2\gamma^{\ga-1}(r_n)}{2+(\ga-1)s\gamma^{\ga-1}(r_n)}\right)^{\frac{\ga}{\ga-1}}.
 $$
From this inequality, we derive
 $$\BA {lll}\myint{0}{t}\ln^\ga(\gw_n(s)+1)ds\leq \myint{0}{t}\left(\myfrac{2\gamma^{\ga-1}(r_n)}{2+(\ga-1)s\gamma^{\ga-1}(r_n)}\right)^{\frac{\ga}{\ga-1}}ds\\[4mm]
 \phantom{\myint{0}{t}\ln^\ga(\gw_n(s)+1)ds}
 \leq 2^{\frac{\ga}{\ga-1}}\gamma(r_n)\myint{0}{t\gamma^{\ga-1}(r_n)}\left(2+(\ga-1)\gt\right)^{-\frac{\ga}{\ga-1}}d\gt.
 \EA$$
 Therefore
 \begin{equation}\label{min0}\BA {lll}A_n(t,x)\geq \gamma (r_n)-\myfrac{N(r_n+|x|)^2}{4t}-N\ln(r_n+|x|)-
\myfrac{N}{2}\ln t\\[4mm]\phantom{A_n(t)---------}
-2^{\frac{\ga}{\ga-1}}\gamma(r_n)
 \myint{0}{t\gamma^{\ga-1}(r_n)}\left(2+(\ga-1)\gt\right)^{-\frac{\ga}{\ga-1}}d\gt.
 \EA \end{equation}

 \nind{\it Step 2: The maximal admissible growth}. We claim that
 \begin{equation}\label{min1}
 \displaystyle\liminf_{\abs x\to \infty}|x|^{-\frac{2}{2-\ga}}\ln \tilde g(|x|)>N^{\frac{1}{2-\ga}} \Longrightarrow
 \lim_{n\to \infty}\underline u_{g_n}(t,x)=\Gf_\infty(t)\qquad\forall (t,x)\in Q^\infty_{\BBR^N}.
 \end{equation}
By replacing $\gt\mapsto (2+(\ga-1)\gt)^{-\frac{\ga}{\ga-1}}$ by its maximal value on $(0,t\gamma^{\ga-1}(r_n))$,
$$ 2^{\frac{\ga}{\ga-1}}\gamma(r_n)\myint{0}{t\gamma^{\ga-1}(r_n)}\left(2+(\ga-1)\gt\right)^{-\frac{\ga}{\ga-1}}d\gt\leq
\gamma^\ga (r_n)t.
$$
Then
  \begin{equation}\label{min2}
 A_n(t,x)\geq \gamma (r_n)-\myfrac{N(r_n+|x|)^2}{4t}-N\ln(r_n+|x|)-
\myfrac{N}{2}\ln t-\gamma^\ga (r_n)t:=B_n(t,x),
 \end{equation}
and
$$\prt_tB_n(t,x)=\myfrac{N(r_n+|x|)^2}{4t^2}-\myfrac{N}{2t}-\gamma^\ga (r_n).
$$
 Thus
  \begin{equation}\label{min3}
  \prt_tB_n(t,x)=0\text{ and } t>0\Longleftrightarrow t:=t_n=\myfrac{N(r_n+|x|)^2}{N+\sqrt{N^2+4N(r_n+|x|)^2\gamma^\ga (r_n)}}.
 \end{equation}
 Therefore $ A_n(t_n,x)$ is bounded from below by the maximum of $B_n(t,x)$ which is achieved for $t=t_n$ and
 $$ \BA {ll}B_n(t_n,x)=\gamma (r_n)-N\ln(r_n+|x|)-\myfrac{N+\sqrt{N^2+4N(r_n+|x|)^2\gamma^\ga (r_n)}}{4}\\[4mm]
 \phantom{udfu"a'g}
 -\myfrac{N(r_n+|x|)^2\gamma^\ga(r_n)}{N+\sqrt{N^2+4N(r_n+|x|)^2\gamma^\ga (r_n)}}
 -\myfrac{N}{2}\ln\left(\myfrac{N(r_n+|x|)^2}{N+\sqrt{N^2+4N(r_n+|x|)^2\gamma^\ga (r_n)}}\right).
\EA $$ Since $r_n\to \infty$ as $n\to \infty$ it follows from last
representation that
  \begin{equation}\label{min4}B_n(t_n,x)=r_n\gamma^{\frac{\ga}{2}} (r_n)\left(\myfrac{\gamma^{1-\frac{\ga}{2}}
  (r_n)}{r_n}-N^{\frac{1}{2}}(1+\nu_n(x))\right),
 \end{equation}
 where $\nu_n(x)\to 0$ as $n\to\infty$ uniformly on any compact
 set in $\Bbb R^N$.
Therefore if $g$ satisfies
\begin{equation}\label{min5}
 \displaystyle\liminf_{\abs x\to \infty}|x|^{-\frac{2}{2-\ga}}\ln \tilde g(|x|)>N^{\frac{1}{2-\ga}},
  \end{equation}
then there holds
\begin{equation}\label{min6}
\underset{n\to\infty}{J_n(t_n,x)}\to
\infty\Longrightarrow\lim_{t_n\to 0}\underline
u_{g_n}(t_n,x)=\infty,
  \end{equation}
  uniformly on compact subsets of $\BBR^N$. We fix $m>0$, denote by $\gl_m$ the first eigenvalue of $-\Gd$ in
  $H^1_0(B_m)$, with corresponding eigenfunction $\gf_m$ normalized by $\sup_{B_m}\gf_m=1$ and set, for $\ge>0$,
  $$W_{m,\ge}(t,x)=e^{-(t+\ge)\gl_m}\Gf_\infty(t+\ge)\gf_m(x)\qquad\forall (t,x)\in Q_\infty^{B_m}.$$
  Then
  $$\BA {ll}
  \prt_t W_{m,\ge}-\Gd W_{m,\ge}+W_{m,\ge}\ln^\ga(W_{m,\ge}+1)=W_{m,\ge}
  \left(\ln^\ga(W_{m,\ge}+1)-\ln^\ga(\Phi_\infty(t+\ge)+1)^{\phantom{a^{a^{a}}}}\!\!\!\!\!\!\!\!\right)\\[2mm]
  \phantom{\prt_t W_{m,\ge}-\Gd W_{m,\ge}+W_{m,\ge}\ln^\ga(W_{m,\ge}+1)}
  \leq 0.
  \EA$$
  Since $\underline u_{g_n}$ increases to the prospective minimal solution $\underline u_{\tilde g}$, it follows due to \eqref{min6} that there exists $n_\ge$ such that
  $$\underline u_{\tilde g}(t_{n_\ge},x)\geq \underline u_{g_{n_\ge}}(t_{n_\ge},x)\geq W_{m,\ge}(t_{n_\ge}+\ge,x)
  \qquad\forall x\in {B_m}.
  $$
 Last inequality in virtue of comparison principle implies
  $$\underline u_g(t,x)\geq W_{m,\ge}(t+\ge,x)\qquad\forall (t,x)\in Q_\infty^{B_m}, t\geq t_{n_\ge}.
  $$
  Letting $\ge\to 0$ yields $\underline u_g\geq W_{m,0}$ in $ Q_\infty^{B_m}$. Since $\lim_{m\to\infty}\gf_m(x)=1$,
  uniformly on any compact subset of $\BBR^N$  and  $\lim_{m\to\infty}\gl_m=0$ we derive $\underline u_{\tilde g}\geqslant \Gf_\infty$ and finally
  $\underline u_{ g}\geqslant \Gf_\infty$. This inequality together with \eqref{m13*} leads to $\underline{u}=\Phi_\infty$. \qeda\medskip

  \nind\Remark In the case $\ga=2$, there holds
  \begin{equation}\label{min7}
\myint{0}{t}\ln^2(\gw_n(s)+1)ds\leq 4\gamma(r_n)\myint{0}{t\gamma (r_n)}(2+\gt)^{-2}d\gt\leq t\gamma(r_n).
  \end{equation}
  Therefore $(\ref{min2})$ is replaced by
  \begin{equation}\label{min8}
 A_n(t,x)\geq\gamma(r_n)-t\gamma^2(r_n) -\myfrac{N(r_n+|x|)^2}{4t}-N\ln(r_n+|x|)-
\myfrac{N}{2}\ln t:=B_n(t,x).
 \end{equation}
 A similarly, there exists $t_n>0$ where $t\mapsto B_n(t,x)$ is maximum and in that case
 $$ \BA {ll}B_n(t_n,x)=\gamma (r_n)-N\ln(r_n+|x|)-\myfrac{N+\sqrt{N^2+4N(r_n+|x|)^2\gamma^2 (r_n)}}{4}\\[4mm]
 \phantom{udfu"a'g}
 -\myfrac{N(r_n+|x|)^2\gamma^2(r_n)}{N+\sqrt{N^2+4N(r_n+|x|)^2\gamma^2 (r_n)}}
 -\myfrac{N}{2}\ln\left(\myfrac{N(r_n+|x|)^2}{N+\sqrt{N^2+4N(r_n+|x|)^2\gamma^2 (r_n)}}\right),
\EA $$
 which yields
   \begin{equation}\label{min9}
B_n(t_n,x)=\gamma (r_n)-r_n\gamma (r_n)(N^{\frac{1}{2}}-\nu_n(x)),
 \end{equation}
where $\nu_n(x)\to 0$ as $n\to \infty$ uniformly on any compact
set in $\Bbb R^N$.
 Thus $B_n(t_n,x)\to -\infty$ as $n\to\infty$. A similar type of computation shows that the expression $I_n(t,x)$ defined in
 $(\ref{m17})$ converges to $0$, whatever is the sequence $\{r_n\}$ converging to $\infty$.

\end{document}